\documentclass[12pt,righttag]{article}

\usepackage{latexsym}
\usepackage{amsfonts}
\usepackage{amsmath}

\newcommand\C{{\mathbb C}}
\newcommand\Z{{\mathbb Z}}
\newcommand\Q{{\mathbb Q}}
\newcommand\A{{\mathbb A}}
\renewcommand\k{{\bf k}}
\renewcommand\char{{\rm  char}\,}
\newcommand\spec{{\mathsf{Spec}\,}}
\newcommand\OO{{\mathcal O}}

\def\DC{\operatorname{DC}}
\def\pr{\operatorname{pr}}
\def\DA{\operatorname{DA}}
\def\DE{\operatorname{DE}}
\def\JA{\operatorname{JA}}
\def\JE{\operatorname{JE}}
\def\JC{\operatorname{JC}}

\def\Center{\operatorname{Center}}

\def\Id{\operatorname{Id}}

\def\Mat{\operatorname{Mat}}

\def\End{\operatorname{End}}

\newtheorem{lmm}{Lemma}
\newtheorem{thm}{Theorem}
\newtheorem{prop}{Proposition}
\newtheorem{conj}{Conjecture}
\newtheorem{remark}{Remark}
\newtheorem{dfn}{Definition}

\title{The Jacobian Conjecture is stably equivalent to
 the Dixmier Conjecture }

\begin{document}

\author{Alexei Belov-Kanel  and Maxim Kontsevich}
\maketitle
\section{Introduction}

The {\bf Jacobian Conjecture} $\JC_n$ in dimension $n\ge 1$ asserts that
{\it for any field $\k$ of characteristic zero any polynomial endomorphism
$\phi$ of the $n$-dimensional affine space $\A^n_\k=\spec
\k[x_1,\dots,x_n]$ over $\k$, with Jacobian $1$:  $$\det \left(
\partial \phi^*(x_i)/ \partial x_j\right)_{1\le i,j\le n}=1$$ is an
 automorphism.} Equivalently, one can say that $\phi$ preserves the
 standard top-degree differential form $dx_1\wedge\dots\wedge
dx_n\in\Omega^n(\A^n_\k)$.

The reference due to this well known problem and related questions can
be found in \cite{Essen}, \cite{BCW}.

 By the Lefschetz principle it is sufficient to consider the case  $\k=\C$.
Obviously, $\JC_n$ implies $\JC_m$ if $n>m$. We denote by $\JC_\infty$ the
stable Jacobian conjecture, the conjunction of conjectures $\JC_n$ for
 all finite $n$. The conjecture $\JC_n$ is obviously true  in the case
$n=1$, and it is open for $n\ge 2$.

\medskip
The {\bf Dixmier Conjecture} $\DC_n$ for integer $n\ge 1$ (see
\cite{Dixmier}) asserts that {\it for any field $\k$ of characteristic
zero any endomorphism of the $n$-th Weyl algebra $A_{n,\k}$ over $\k$
is an automorphism.}

Here $A_{n,\k}$ is the associative unital algebra over $\k$ with $2n$
generators $y_1,\dots,y_{2n}$ and relations

 $$[y_i,y_j]=\omega_{ij}\,,$$
 where $(\omega_{ij})_{1\le i,j\le 2n}$ is the following standard
 $2n\times 2n$ skew-symmetric matrix:
 $$\omega_{ij}=\delta_{i,j+n}-\delta_{i+n,j}\,\,.$$ 
 The algebra $A_{n,\k}$
 coincides with the algebra $D(\A^n_\k)$ of  polynomial differential
 operators on $\A^n_\k$.  For any $i,\,1\le i\le n$ element $y_i$ acts
 as the multiplication operator by the variable $x_i$, and element
 $y_{n+i}$ acts by the differentiation $\partial/\partial x_i$.  Again, it
 is sufficient to consider the case  $\k=\C$. The conjecture $\DC_n$ implies
 $\DC_m$ for $n>m$, and we can consider the stable Dixmier conjecture
 $\DC_\infty$. The conjecture $\DC_n$ is open for any $n\ge 1$.

It is well-known that $\DC_n$ implies $\JC_n$ (in particular $\DC_\infty$
implies $\JC_\infty$) (see \cite{Essen}, \cite{BCW}).  The argument is
very easy. Let $\phi:\A^n_\k\to \A^n_\k$ be a counterexample to $\JC_n$.
Then $\phi$ is a non-invertible \'etale map, and it induces a pullback
homomorphism $\phi^*_{diff}$ of the algebra of differential
operators on $\A^n_\k$. The endomorphism $\phi^*_{diff}$ of the Weyl
algebra preserves the degree of differential operators. Restricting
$\phi^*_{diff}$ to zero order differential operators, we obtain
the usual pullback $\phi^*$ of functions on $\A^n_\k$.  By our
assertion it is not surjective, hence we obtain a counterexample to
$\DC_n$.

   Our result  is an opposite implication. Namely, we prove the following

   \begin{thm}
   Conjecture $\JC_{2n}$ implies $\DC_n$.
      \end{thm}

 In particular, we obtain  that the stable conjectures $\JC_\infty$ and
$\DC_\infty$ are equivalent.

 \begin{remark}
 A.~van den Essen (\cite{Essen}, Theorem 10,4.2) proved a
  weaker result: the conjecture $\JC_{2n}$ implies the invertibility of any
  endomorphism of $A_{n,\k}=D(\A^n_\k)$ preserving the filtration by
  the degrees of differential operators.
  \end{remark}

For the convenience of the reader, and in order to make the text
  self-contained, we include in the paper proofs of several known
results scattered in the literature.  During the preparation of this
paper we have learned from K.~Adjamagbo about the preprint \cite{Tsuch}
where two key results concerning the Weyl algebra in finite characteristic
were established (Propositions 2 and 4 from Section 4 in the present
paper), see also a very recent preprint \cite{ACE}.

\begin{remark} The present paper is written in the 
standard language of algebraic geometry.
It is possible (and reasonable for some minds) to use the
model-theoretic language of non-standart analysis, instead of
sch\'eme-theoretic considerations. In particular, in the proofs of
 several results of  
our paper one can
use the reduction modulo an infinitely large prime.    \end{remark}

\begin{remark} After this paper was written, we were told by Ken Goodearl about a paper
"Endomorphsims of Weyl algebra and p-curvatures"
(Osaka Journal of Mathematics Volume 42, Number 2 (June 2005)) by 
Yoshifumi Tsuchimoto which contains the proof of our main result.
 The proofs by  Tsuchimoto and in the present paper are different (although there are many similarities),
   hence we think that it is 
 reasonable to keep our paper on archive.

\end{remark}

{\bf Acknowledgments}: we are grateful to Kossivi Adjamagbo, Jean--Yves
Charbonnel, Ofer Gabber and Leonid Makar--Limanov for useful
discussions and comments.

\section{A reformulation of the Jacobian conjecture}

For  given integers $n\ge 2,d\ge 1$ we denote by $\JE_{n,d}$ an affine
scheme of finite type over $\Z$ representing the following functor.
For any commutative ring $R$  the set $\JE_{n,d}(R)$ is the set of
endomorphisms $f$ of $R$-algebra $R[x_1,\dots,x_n]$ such that

\begin{itemize}
\item $\det\left( \partial f(x_i)/ \partial x_j\right)_{1\le i,j\le
n}=1\in R[x_1,\dots,x_n],$

\item $ \deg (f(x_i))\le d\,\,\,\,\,\,\forall i,\,1\le i\le n\,.$
\end{itemize}

We say that $f$ as above is an {\it endomorphism of  degree $\le d$}
(and with Jacobian $1$).  The ring of functions $\OO(\JE_{n,d})$ is
finitely generated, its  generators are coefficients $c_{i,\alpha}$
 which appear in the universal endomorphism over $\OO(\JE_{n,d})$:
$$f_{univ}(x_i)=\sum_{\alpha:\,|\alpha|\le d} c_{i,\alpha} x^\alpha$$

Here $\alpha=(\alpha_1,\dots,\alpha_n)\in\Z_{\ge 0}^n$ is a
multi-index, $x^\alpha:=\prod_{i=1}^n x_i^{\alpha_i},\,\,
 |\alpha|:=\sum_{i=1}^n \alpha_i\,.$

Similarly, for $n\ge 2,d\ge 1,d'\ge 1$ we denote by $\JA_{n,d,d'}$ an
 affine scheme of finite type over $\Z$ parameterizing pairs of
endomorphisms $(f,f')$ of $n$-dimensional affine space, with Jacobian
$1$, of degrees $\le d$ and $\le d'$ respectively, mutually inverse to
each other:  $f\circ f'=f'\circ f=\Id_{\A^n}$.

We have an obvious forgetting map $\pr_{n,d,d'}^{(J)}:
\JA_{n,d,d'}\to \JE_{n,d},\,\,(f,f')\mapsto f$ which is an immersion
(i.e. $\JA_{n,d,d'}$ is identified with a locally closed subscheme of
$\JE_{n,d}$).  The Jacobian conjecture $\JC_n$ means that for any $d\ge
1$
$$
\JE_{n,d}\times\spec \Q=\bigcup\limits_{d'\ge 1}
\pr_{n,d,d'}^{(J)}(\JA_{n,d,d'}\times \spec \Q)\,\,.
$$

For given $n,d$ the set $X_{d'}:=(\JA_{n,d,d'}\times \spec \Q)\subset
\JE_{n,d}\times \spec\Q$ is a constructible set.  Therefore, we get an
infinite growing chain of  constructible subsets $X_1\subset
X_2\subset\dots$ of the scheme of finite type $\JE_{n,d}\times \spec \Q$
over $\Q$.

Let us assume $\JC_n$ and fix an integer $d\ge 1$. Then $\cup_{d'\ge 1}
X_{d'}= X$ where $X:=\JE_{n,d}\times \spec\Q$.  Then it follows from the
standard  properties of constructible sets (see \cite{EGA}, Corrollaire
1.9.8, Chapitre IV) that there exists an integer $d'$ such that $X_{d'}=X$.
Alternatively, one can use a result of O.~Gabber (see \cite{BCW},
Theorem 1.2) which says that for an automorphism $f$ of
$\k[x_1,\dots,x_n]$ of degree $\le d$ in the above sense ($\k$ is a
field of any characteristic),  the inverse map has the degree $\le
d^{n-1}$. Hence one can a priori set $d'=d^{n-1}$. Anyhow,  the
Jacobian conjecture can be rephrased as the equality
$\JA_{n,d,d'}(\C)=\JE_{n,d}(\C)$.

   The following statement is obvious.

\begin{lmm} \label{lmm:0-finite}
Let $\phi:A\to B$ be an immersion of schemes of finite type over
$\Z$. Then $\phi$ induces a bijection between $A(\C)$ and $B(\C)$ if
and only if there exists a finite set of primes $S$ such that $\phi$
induces a bijection between $A(\k)$ and $B(\k)$ for any  field $\k$
with $\char \k\notin S\cup \{0\}$.
\end{lmm}

We apply it to the projection $\pr_{n,d,d'}^{(J)}$. The
conclusion is that the Jacobian conjecture $\JC_n$ is equivalent to the
following

\begin{conj}{($\JC_n$ in finite characteristic)}\label{conj:JCf}
For any $d\ge 1$ there exists $d'\ge 1$ and a finite set of primes $S$
such that for any field $\k$ with $\char\k\notin S\cup\{0\}$ and any
polynomial map $\phi:\A_\k^n\to\A_\k^n$ of degree $\le d$ with Jacobian
$1$, the inverse map exists and has degree $\le d'$.
\end{conj}

The equivalence of $\JC_n$ and the above conjecture in finite
characteristic was first established by K.~Adjamagbo in
\cite{Adjamagbo}.

\section{More about Weyl algebras}

\subsection{Weyl algebras over an arbitrary base}\label{sec:WA}

One can define an algebra $A_{n,R}$ for arbitrary commutative ring $R$
exactly in the same way as for fields of characteristic zero. This
algebra  is free as a $R$-module. It has a canonical basis consisting of
monomials $(y_1^{\alpha_1}\dots
y_{2n}^{\alpha_{2n}})_{(\alpha_1,\dots,\alpha_{2n})\in\Z_{\ge
0}^{2n}}$.  Although for any $R$ the algebra $A_{n,R}$ maps to the algebra
$D(\A_{R}^n)$ of differential operators acting on $R[x_1,\dots,x_n]$,
these two algebras are not isomorphic in general. For example, if $R$
is an algebra over $\Z/p\Z$ for some prime $p$, then the operator $(d/dx_1)^p$
is zero.

We say that an endomorphism $f$ of the algebra $A_{n,R}$ {\it has degree}
$\le d$ if the image $f(y_i)$ of any generator $y_i\in
A_{n,R},\,\,\,i=1,\dots,2n$ is a linear combination with coefficients
in $R$  of the monomials of degree $\le d$.  In a manner completely
parallel to the previous section, we can define schemes of finite type
$\DE_{n,d}$, $\DA_{n,d,d'}$, and the projection $\pr^{(D)}_{n,d,d'}$.
Also, we can make a reformulation of the Dixmier conjecture in the same
way as for the Jacobian conjecture.

 \subsection{The Weyl algebra in finite characteristic as an Azumaya
algebra}

It is a classical fact that in finite characteristic  the algebra $A_{n,R}$
has a big center, and it is moreover an Azumaya algebra of its center
 (see \cite{Revoy}).

We will use the following slightly non-standard definition of an Azumaya algebra
(see e.g. \cite{Milne},  Proposition 2.1, Chapter IV):

 \begin{dfn}
For a commutative ring $R$ an Azumaya algebra  over $R$ of rank $N\ge
1$ is an associative unital algebra $A$ over $R$ which is a finitely
generated $R$-module and such that there exists a finitely generated
faithfully flat extension $R'\subset R$ of $R$ such that the pullback
algebra $A':=A\otimes_R R'$ is isomorphic to the matrix algebra
$$\Mat(N\otimes N,R')=\Mat(N\times N,\Z)\otimes R'$$ as an
algebra over $R'$.
 \end{dfn}

It follows by descent that the center of an Azumaya algebra over $R$
coincides with $R$.  Also, an Azumaya algebra $A$ considered as a 
$R$-module is a finitely generated projective module, in other words, a
vector bundle over $\spec R$. This bundle has rank $N^2$, its fibers
are associative algebras, and the  fiber over any point of $\spec R$
over an algebraically closed field $\k$ is isomorphic to the matrix
algebra $\Mat(N\times N,\k)$.

\begin{prop}
For any commutative algebra $R$  over $\Z/p \Z$ where $p$
is a prime, the algebra $A_{n,R}$ is an Azumaya algebra of rank $p^n$ over
$R[x_1,\dots,x_{2n}]$.  The central element of $A_{n,R}$ corresponding
to variable $x_i$ is $y_i^p$.
\end{prop}

{\bf Proof}: Let us introduce a faithfully flat extension
$R':=R[\xi_1,\dots,\xi_{2n}]$ of $R[x_1,\dots,x_{2n}]$, where the
inclusion of $R[x_1,\dots,x_{2n}]$ into $R'$ is given by
 $$
 x_i\mapsto
 \xi^p,\,\,\,\,i\in\{1,\dots, 2n\}\,.
 $$

We claim that the algebra over $R'$
$$A':=A_{n,R}\otimes_{R[x_1,\dots,x_{2n}]} R[\xi_1,\dots,\xi_{2n}]$$
is isomorphic to the matrix algebra of rank $p^n$ over $R'$.
Namely, the algebra $A'$ considered as an algebra over
 $R'=R[\xi_1,\dots,\xi_{2n}]$, has generators $y_i,\,i\in
\{1,\dots,2n\}$ and defining relations
$$ [y_i,y_j]=\omega_{ij},\,\,\, y_i^p=\xi_i^p \,\,.$$
Introduce a new set of generators $y'_i\in A',\,\,\,i\in
\{1,\dots,2n\}$ by the formula
$$y_i':=y_i-\xi_i\,\,.$$
These generators have  defining relations
$$[y_i',y_j']=\omega_{ij},\,\,\, (y_i')^p=0\,\,. $$

Hence, we see that the algebra $A'$ over $R'$ is isomorphic to the tensor
product over $\Z/p\Z$ of $R'$ and a finite-dimensional algebra over
$\Z/p \Z$ given by the generators $(y_i')_{1\le i\le 2n}$ and the 
relations as
above. The last algebra is the tensor product of $n$ copies of its
version in the case $n=1$.  The statement of the proposition now
follows from the following

\begin{lmm}
For any prime number $p$ the 
algebra $A$ over $\Z/p\Z$ with two generators $y_1,y_2$ and relations
 $$[y_1,y_2]=1,\,y_1^p=y_2^p=0$$
is isomorphic to $\Mat(p\times p, \Z/p\Z)$.
\end{lmm}

{\bf Proof}: Consider the finite ring $B:=\Z/p\Z[x]/(x^p)=\Z[x]/(x^p,p)$.
It is isomorphic to $(\Z/p\Z)^p$ as an abelian group (and as
$\Z/p\Z$-module).  Differential operators $Y_1,Y_2$ acting on $B$ and
given by the formulas $$Y_1(f)=df/dx,\,\,\,Y_2 f= xf,\,\,f\in B$$ are
well-defined, and satisfy the relations $[Y_1,Y_2]=1,\,Y_1^p=Y_2^p=0$.
Hence we obtain a homomorphism $A\to \End_{\Z/p\Z-mod}(B)$.  A direct
calculation shows that this is an isomorphism.  $\Box$

    \section{The proof of the implication $\JC_{2n}\to \DC_n$}

Let us assume that the conjecture $\JC_{2n}$ (phrased in the form of
Conjecture~\ref{conj:JCf}) is true, our goal is to prove $\DC_n$.

Let $f:A_{n,\C}\to A_{n,\C}$ be an endomorphism of degree $\le d$. We
have to prove that $f$ is invertible.

Denote by $R$ the subring of $\C$ generated by the coefficients of elements
$f(y_i)\in A_{n,\C},\,,i\in\{1,\dots,2n\}$ in the standard basis of
$A_{n,\C}$.  The ring $R$ is a finitely generated integral domain.
Moreover, we may assume  that for any prime $p$ the 
ring $R/p R$ is either
zero or an integral domain, in particular it has no non-zero
nilpotents. In order to achieve this property it is enough to extend
$R$ by adding inverses to finitely many primes.

For any prime $p$ the endomorphism $f$ induces an endomorphism
$$f_p:A_{n,R/p R}\to A_{n,R/p R}$$
of an Azumaya algebra of rank $p^n$ over
$$C_p:=R/p R[x_1,\dots, x_{2n}]=\Center(A_{n,R/pR})\,\,.$$

The following result was proved first by Y.~Tsuchimoto \cite{Tsuch},
it follows also from a more general recent result from \cite{ACE}.

\begin{prop}
The endomorphism $f_p$ maps $C_p$ to itself.
\end{prop}

{\bf Proof}: Denote by $\k$  an algebraically closed field of
characteristic $p$. For any $\k$-point $v$ of $\spec C_p$ the fiber
$A_v$ is an algebra over $\k$ isomorphic to $\Mat(p^n\times p^n,\k)$.

\begin{lmm}
An element $a\in A_{n,R/p R}$ belongs to $C_p$ if and only if
for any $\k$-point $v$ of $\spec C_p$ where $\k$ is an algebraically
closed field, the image of $a$ in $A_v$ is central, i.e. it is a scalar
matrix.
\end{lmm}

{\bf Proof}: One direction is obvious, i.e. if $a$ is central than its
image in $A_v$ is central. Conversely, if $a\in A_{n,R/p R}$ is not
central then there exists $b\in  A_{n,R/p R}$ such that $[a,b]\ne 0$.
For any non-zero section $s$ of the vector bundle $A_{n,R/pR}/C_p$ there
exists a $\k$-point at which this section does not vanish, because
the algebra $C_p=R/p R[x_1,\dots, x_{2n}]$ has no non-zero nilpotents by
our assumption that $R/p R$ has no non-zero nilpotents. We apply this
argument to the section $s=[a,b]$ and conclude that the image of $a$ in
$A_v$ is not central for some $v$.  $\Box$

Let $a\in C_p\subset A_{n,R/pR}$ be a central element. We want to prove
that $f(a)$ is central. Assume the opposite. Then by the above lemma there
exists a homomorphism $\rho:A_{n,R/pR}\to \Mat(p^n\times p^n,\k)$
such that $\rho(f(a))$ is not a scalar matrix.  Let us denote by
$V_0\simeq \k^{p^n}$ the module over $A_{n,R/pR}$ associated to the
homomorphism $f\circ \rho$.  Our assumption mean that $V_0$ considered
as a module over $\k\otimes C_p$  is {\it not} isomorphic to the sum of
$p^n$ copies of the simple module $M_v\simeq \k$ associated with any
$\k$-point $v$ of $\spec C_p$. The support of the module $V_0$ is a
non-empty finite subscheme of $\spec C_p$ defined over $\k$, hence
there exists a $\k$-point $v$ in it support. Moreover, for any such
point $v$ the tensor product $V:=V_0\otimes_{C_v\otimes \k} M_v$ is a
vector space over $\k$ such that $0<\dim V<\dim V_0$.  The algebra
$A_{n,R/pR}$ maps to the algebra of endomorphisms of $C_p$-module
$V_0$, hence it maps to the algebra of $\k$-linear endomorphisms of
$V$. In this representation of $A_{n,R/pR}$ the center $C_p$ acts by
scalars, by the nature of the construction.

Therefore, we obtain a homomorphism $C_p\to \k$, i.e. 
a $\k$-point $v$ of $\spec C_p$,
and a homomorphism of $\k$-algebras
$$
A_v\simeq \Mat(p^n\times p^n, \k)\to
\Mat(M\times M,\k),\,\,\,0<M<p^n\,,
$$
here $M:=\dim V$.
This is impossible because $\Mat(p^n\times p^n,\k)$ is simple and
$0<\dim_\k \Mat(M\times M,\k)< \dim_\k \Mat(p^n\times p^n,\k)$.
We obtain a contradiction. The Proposition is proven. $\Box$

Denote by $f^{centr}_p$ the endomorphism of
$C_p$ induced by $f$.
Our next goal is to prove that $f^{centr}_p$ {\bf  preserves
certain $R/pR$-linear Poisson bracket} on $C_p$.

Namely, we define an operation
$\{\,,\,\}:C_p\otimes_{R/pR} C_p\to C_p$
by the formula
$$
\{a,b\}=\frac{[\tilde{a},\tilde{b}]}{p}\pmod
{pA_{n,R}}\in A_{n,R/pR}=A_{n,R}/pA_{n,R}
$$
where $\tilde{a},\tilde{b}\in A_{n,R}$ are arbitrary lifts of the elements
$a,b\in C_p\subset A_{n,R/pR}$.  First of all, it is easy to see
that the commutator $[\tilde{a},\tilde{b}]$ vanishes modulo $p$, hence
the division by $p$ makes sense. It is uniquely defined because $R$ and
hence $A_{n,R}$ both have no torsion.  A straightforward check shows
that $\{a,b\}$   defined as above does not depend on the choice of
the lifts $\tilde{a},\tilde{b}$, and it belongs to the center $C_p$.
Moreover, the  commutator $\{\,,\,\}$ on $C_p$ is a $R/pR$-linear,
skew-symmetric operation satisfying the Jacobi identity (hence $C_p$
becomes a Lie algebra), and for any $a\in C_p$ the operator
$\{a,\cdot\}:C_p\to C_p$ is a $R/pR$-linear derivation of $C_p$, i.e.
the bracket satisfies the Leibniz rule
 $$\{a,bb'\}=\{a,b\}b'+\{a,b'\}b\,\,.$$

\begin{lmm}
The above defined canonical
Poisson bracket on $C_p\simeq R/p[x_1,\dots,x_{2n}]$
is given by the standard formula
$$
\{a,b\}=\sum_{i=1}^n \left(
\frac{\partial a}{\partial x_i}\frac{\partial b}{\partial x_{n+i}}
 - \frac{\partial b}{\partial x_i}\frac{\partial a}{\partial
x_{n+i}}\right)\,\,.
$$
\end{lmm}

{\bf Proof}: By the Leibniz rule it follows that it suffices to calculate
the bracket $\{x_i,x_j\}$ for any two generators of $C_p$.  The
calculation reduces to the case $n=1$.  It is convenient to calculate
first the commutator in the algebra $A_{1,\Z}$ and then make the 
reduction modulo $p$:
$$
\frac{1}{p}[(d/dx)^p,x^p]=\frac{1}{p}\sum_{i=0}^{p-1}
\frac{(p!)^2}{(i!)^2 (p-i)!} x^i (d/dx)^i=
- 1\pmod {p}
$$
Then the statement of the lemma follows immediately. $\Box$

The next lemma follows directly from the  definition of the bracket:

\begin{lmm}
The homomorphism $f^{centr}_p:C_p\to C_p$
preserves the canonical Poisson bracket.
\end{lmm}

It is well-known in symplectic geometry that a non-degenerate Poisson
structure on a $C^{\infty}$ manifold $X$ is essentially the same as a
symplectic structure, i.e. a non-degenerate closed $2$-form.  The same
 is true in the algebraic context, in characteristic $>2$. Namely, a 
Poisson bracket gives a section
$$\alpha\in \Gamma(\spec C_p,\wedge^2 T_{\spec C_p/\spec R/pR})$$
of the wedge square of the tangent bundle, defined by the formula
$$
\{f,g\}=\langle df\wedge dg, \alpha\rangle \in C_p,\,\,\,\forall
f,g\in C_p\,\,.
$$
This section can be interpreted as an operator from
the cotangent bundle to the tangent bundle.  This operator is
invertible in our case, the inverse operator can be interpreted as a
$2$-form
$$
\omega:=\alpha^{-1}=\frac{1}{2}\sum_{1\le i,j\le 2n}
\omega_{ij} dx_i\wedge dx_j=\sum_{i=1}^n dx_i\wedge dx_{n+i}\,\,.
$$

\begin{lmm}
For $p>n$ the endomorphism $f^{centr}_p$ of $C_p= R/p R[x_1,\dots, x_{2n}]$
preserves the top-degree form
$dx_1\wedge\dots \wedge dx_{2n}\in \Omega^{2n}(B/(R/pR))$.
\end{lmm}

{\bf Proof}: It follows from the previous lemma that $f_{centr}$ preserves
the symplectic $2$-form
$\omega$.
The volume form from above is equal
to $\pm\omega^n/n!$ for $p>n$. $\Box$

The next result implies that the degree of $f^{centr}_p$ is $\le d$.

\begin{prop}\label{prop:degree}
For any field $\k$ of characteristic $p$ and any $\k$-point $v$ of
$\spec R$, the degree of $f_v$ (as an endomorphism of the Weyl algebra
$A_{n,\k}$) is equal to the degree of $f^{centr}_v$ (as an endomorphism
of the polynomial algebra) where $f^{centr}_v$ is the endomorphism of
$\Center(A_{n,\k})\simeq\k[x_1,\dots,x_{2n}]$ induced from
$f_p^{centr}$.
\end{prop}

{\bf Proof}: The degree of $f^{centr}_v$ is defined as the maximum
over $i\in\{1,\dots,2n\}$ of the degrees of polynomials
$f^{centr}_v(x_i)$.  The degree of endomorphism $f_v$ is defined as
the maximum over $i\in\{1,\dots,2n\}$ of the degrees (in the sense of
Bernstein filtration, by the degree of monomials in the standard basis
of $A_{n,\k}$) of elements $f_v(y_i)$.  We claim that for each index
$i$ both degrees coincide with each other.  The reason is the
following.  Let $d_i$ be the degree of $f_v(y_i)$.  We claim that the 
degree of $f_v(y_i^p)=(f_v(y_i))^p$ considered as an element of
$A_{n,\k}$, is equal to $p d_i$.  It follows from the following
\begin{lmm}\label{lmm:degree} The degree is an additive character of the
multiplicative monoid of non-zero elements in $A_{n,\k}$.  \end{lmm}
{\bf Proof}: It follows immediately from the consideration of Bernstein
filtration on $A_{n,\k}$ and the remark that the product of non-zero
homogeneous polynomials is a non-zero polynomial.  $\Box$

The degree of $f_v(y_i^p)$ considered as an element of
$\Center(A_{n,\k})$ is $1/p$ times its degree in $A_{n,\k}$, i.e.
it is equal to $d_i$.  Proposition \ref{prop:degree} is proven. $\Box$

Now we can use finally our main assumption that the Jacobian conjecture
$\JC_{2n}$ holds.  Namely, by its reformulation (in form of Conjecture
\ref{conj:JCf}), we conclude that there exists an integer $d'\ge 1$ and
a finite set of primes $S$ (the union of the set of excluded primes for
$\JC_{2n}$ in form of Conjecture \ref{conj:JCf}, and the set of primes
$\le n$), such that for any algebraically closed field $\k$ such that
$p=\char (\k)\notin S\cup\{0\}$ and any $\k$-point $v$ of $\spec R$,
the pullback $f^{centr}_v$ to $v$ of $f^{centr}_p$ is
invertible and the inverse endomorphism of $\k[x_1,\dots,x_{2n}]$ has
the degree $\le d'$.

The following result is a particular case of a more general statement
proven in \cite{ACE}, and also follows from \cite{Tsuch}.

\begin{prop}
For any $v$ as above the endomorphism $f_v$ of $A_{n,\k}$ is invertible.
\end{prop}

{\bf Proof}:  We may assume that $\k$ is algebraically closed.
The endomorphism $f_v$ of Azumaya algebra $A_{n,\k}$ preserves the
center and is invertible on the center. Thus, it gives a $C_v$-linear
homomorphism $g_v$ from  one Azumaya algebra of rank
$p^n$ over $C_v:=\k[x_1,\dots,x_{2n}]$ (here we mean the
algebra $A_{n,\k}$), to another Azumaya algebra of rank
$p^n$ (the pullback of $A_{n,\k}$ by $f_v^{centr}$).
Any such a homomorphism  restricts to an isomorphism after the 
reduction to any
$\k$-point of $C_v$, because any homomorphim of associative $\k$-algebras
$$\Mat(N\times N,\k)\to \Mat(N\times N,\k),\,\,N:=p^n$$
is an isomorphism. Therefore, $g_v$ is an isomorphism
of vector bundles. $\Box$

Finally, the degree of the inverse to $f_v$ is $\le d'$, as follows
directly from Proposition \ref{prop:degree}.  The conclusion is that
for any point $v$ of $\spec R$ over an field $\k$ of finite
characteristic $p\notin S$, the corresponding point of the scheme of finite
type $\DE_{n,d}$ (see Section \ref{sec:WA} for the notation) belongs to
the constructible set $\DA_{n,d,d'}$.  This implies (see Lemma
\ref{lmm:0-finite}) that $f$ is invertible after the localization to
zero characteristic, and the inverse endomorphism has degree $\le d'$.
Theorem 1 is proven. $\Box$

\begin{remark} It is interesting that Poisson brackets appear in
another situation related to polynomial automorphisms. The Poisson
algebra structure was used by I.~Shestakov and U.~Umirbaev in their
proof that the Nagata automorphism is wild (see \cite{ShUm}).  \end{remark}

\vspace{5mm}

{\bf Addresses:}

A.B.-K.: Institute of Mathematics, Hebrew University, Givat Ram,
Jerusalem 91904, Israel.

{kanel@mccme.ru} \\

M.K.: IHES, 35 route de Chartres, Bures-sur-Yvette 91440, France

{maxim@ihes.fr}


\begin{thebibliography}{99}


\bibitem{Adjamagbo}
K.~Adjamagbo,  {\it  On separable algebras over a U.F.D. and the Jacobian
conjecture in any characteristic. Automorphisms of affine spaces} (Curacao,
1994), 89--103,
Kluwer Acad. Publ., Dordrecht, 1995.
(prepublication 91, 1996, institute de matematiques de Jussieu.)


\bibitem{ACE} K.~Adjamagbo,  J.~Y.~Charbonnel, A.~van den Essen,
{\it On ring homomorphisms of Azumaya algebras}, e-print math/0509188.

\bibitem{BCW}
 H.~Bass,  E.~H.~Connell, D.~Wright,
{\it The Jacobian conjecture: reduction of degree and formal expansion
of the inverse.}  Bull. Amer. Math. Soc. (N.S.) {\bf 7 }(1982), no. 2,
287--330.

\bibitem{Dixmier}
J.~Dixmier, {\it  Sur les algebres de Weyl},
Bull. Soc. Math. France {\bf 96} (1968), 209--242.

\bibitem{Essen}
A.~van den Essen,{\it Polynomial automorphisms and the Jacobian
conjecture,} Progress in Mathematics, 190. Birkhauser Verlag, Basel,
2000.

\bibitem{EGA} A.~Grothendieck, J.~Dieudonne , {\it
El\'ements de G\'eometrie Alg\`ebrique. IV.
Etude locale des sch\'emas et des morphismes de sch\'emas}, Publ. Math.
Inst. Hautes Etudes Sci.
{\bf 20} (1964).

\bibitem{Milne} J.~S.~Milne, {\it \'Etale Cohomology},
Princeton Mathematical Series, 33, Princeton University Press, 1980.

\bibitem{Revoy} P.~Revoy, {\it Alg\`ebres de Weyl en charact\'eristique p},
Compt. Rend. Acad. Sci. Paris,
S\'er. A-B {\bf 276} (1973), A, 225--228.

\bibitem{ShUm}
I.~P.~Shestakov; U.~U.~Umirbaev,
{\it The tame and the wild automorphisms of polynomial rings in three
variables},  J. Amer. Math. Soc. {\bf 17} (2004), no. 1, 197--227.

\bibitem{Tsuch}
Y.~Tsuchimoto, {\it Preliminaries on Dixmier conjecture},
Mem. Fac. Sci. Kochi Univ. Ser. A Math. {\bf 24} (2003), 43--59.

\end{thebibliography}
    \end{document}